\pdfoutput=1
\documentclass[11pt]{article}

\usepackage[a4paper,margin=1.1in]{geometry}
\usepackage[T1]{fontenc}
\usepackage[utf8]{inputenc}
\usepackage{lmodern}
\usepackage{microtype}
\usepackage[numbers,sort&compress]{natbib}
\usepackage{amsmath,amssymb,amsthm,mathtools}
\usepackage{url}
\usepackage[colorlinks=true,linkcolor=blue,citecolor=blue,urlcolor=blue]{hyperref}
\hypersetup{
  pdftitle={Real Matrices as Linear Combinations of Three Orthogonal Matrices},
  pdfauthor={Zhekai Pang}
}

\theoremstyle{plain}
\newtheorem{lemma}{Lemma}
\newtheorem{theorem}{Theorem}
\newtheorem{corollary}{Corollary}

\newcommand{\R}{\mathbb{R}}

\DeclareMathOperator{\diag}{diag}

\title{Real Matrices as Linear Combinations of Three Orthogonal Matrices}
\author{Zhekai Pang\thanks{Universitat Pompeu Fabra, Barcelona, Spain. Email: \texttt{zhekai.pang@upf.edu}.}}
\date{}

\begin{document}

\maketitle

\begin{abstract}
Li and Poon proved that every real square matrix is a real linear combination of four real orthogonal matrices. The resulting question, recorded by Zhan, asks whether four terms are necessary. We prove that three real orthogonal matrices always suffice.
\end{abstract}



\section{Introduction}

Additive decomposition problems ask how efficiently matrices can be expressed as linear combinations of elements from a prescribed structured class. For unitary and orthogonal matrices, this question is particularly natural: each summand preserves the underlying inner-product geometry, whereas a general matrix may distort lengths and angles, and may have arbitrary singular values. The problem is therefore to determine how many such rigid building blocks are needed to represent arbitrary matrices.

The complex and real problems have different behavior. In the complex case, every square complex matrix is a linear combination of two unitary matrices; this classical result is also discussed by Wu \cite{Wu1994} and Zhan \cite{Zhan2001}. In the real case, Zhan proved that every real matrix of order $n$ is a real linear combination of $n$ real orthogonal matrices \cite{Zhan2001}. Li and Poon later proved that every real square matrix is a real linear combination of four real orthogonal matrices \cite{LiPoon2002}. This led to the question, recorded by Zhan, of whether four terms are least possible \cite[Question~22]{Zhan2007}. Li, Zhang, and Li \cite{LiZhangLi2022} later proved that at least three terms are needed in every dimension at least three.

We prove that three terms always suffice. Combining this upper bound with the lower bound of Li, Zhang, and Li \cite{LiZhangLi2022} gives the exact answer to Zhan's question: the minimum number is one in dimension one, two in dimension two, and three in every dimension at least three.

\section{Preliminaries}

Let $M_n(\R)$ be the real vector space of $n\times n$ real matrices. Let
\[
O(n)=\{Q\in M_n(\R):Q^TQ=I_n\}
\]
be the real orthogonal group. For a positive integer $k$, define
\[
\mathcal A_n(k)=\{\alpha_1Q_1+\cdots+\alpha_kQ_k:
\alpha_i\in\R,\ Q_i\in O(n)\ \text{for } i=1,\ldots,k\}.
\]
Coefficients equal to zero are allowed. Thus $\mathcal A_n(k)$ is the set of matrices that are real linear combinations of at most $k$ real orthogonal matrices.

For $A\in M_n(\R)$, write $s_1(A)\geq\cdots\geq s_n(A)
\geq0$ for the singular values of $A$, counted with multiplicity.

The following elementary sufficient condition is a special case of the two-term characterization of Li, Zhang, and Li \cite{LiZhangLi2022}.

\begin{lemma}\label{lem:two-term}
Let $n\geq3$, and let $X\in M_n(\R)$.

If $n=2m+1$ with $m\geq1$, and the singular values of $X$, counted with multiplicity, can be listed as
\[
R,t_1,t_1,t_2,t_2,\ldots,t_m,t_m,
\]
where $R\geq t_1\geq\cdots\geq t_m\geq0$, then $X\in\mathcal A_n(2)$.

If $n=2m$ with $m\geq2$, and the singular values of $X$, counted with multiplicity, can be listed as
\[
R,t_1,t_1,t_2,t_2,\ldots,t_{m-1},t_{m-1},r,
\]
where $R\geq t_1\geq\cdots\geq t_{m-1}\geq r\geq0$, then $X\in\mathcal A_n(2)$.

Moreover, the construction in Appendix A writes $X=\alpha U+\beta V$ with $U,V\in O(n)$, where
\[
\alpha=\frac{R+r_X}{2},\qquad \beta=\frac{R-r_X}{2}.
\]
Here $r_X=s_n(X)$ is the smallest singular value of $X$, so $r_X=t_m$ in the odd-dimensional case and $r_X=r$ in the even-dimensional case.
\end{lemma}

The proof is given in Appendix A.

The next lemma is a standard inverse form of a rank-one update result; see Golub \cite[Section~5]{Golub1973} and Bunch, Nielsen, and Sorensen \cite{BunchNielsenSorensen1978}.

\begin{lemma}\label{lem:rank-one}
Let $\ell\geq2$, let $x_1>\cdots>x_\ell$ be real numbers, and let $\sigma>0$. Let $y_1,\ldots,y_{\ell-1}$ satisfy
\[
x_1>y_1>x_2>y_2>\cdots>x_{\ell-1}>y_{\ell-1}>x_\ell.
\]
Let $L=x_1+\cdots+x_\ell+\sigma-y_1-\cdots-y_{\ell-1}$, and assume that $L>x_1$.
Then there is a unit vector $w\in\R^\ell$ such that the eigenvalues of $\diag(x_1,\ldots,x_\ell)+\sigma ww^T$, counting multiplicities, are exactly $L,y_1,\ldots,y_{\ell-1}$.
\end{lemma}

The proof is given in Appendix A.

\section{Main result}

\begin{lemma}\label{lem:strict-diagonal}
Let $n\geq3$, and let
\[
D=\diag(d_1,\ldots,d_n),\qquad d_1>d_2>\cdots>d_n\geq0.
\]
Let $c>d_1$. Then there is a unit vector $u\in\R^n$ such that, with $H=I_n-2uu^T$ and $S=D-cH$, the matrix $S$ satisfies one of the two singular-value hypotheses in Lemma \ref{lem:two-term}. Consequently, $D\in\mathcal A_n(3)$.
\end{lemma}

\begin{proof}
Let $a_i=d_i-c$ for $i=1,\ldots,n$. Since $c>d_1$, we have $0>a_1>a_2>\cdots>a_n$. For a unit vector $u\in\R^n$, let $H=I_n-2uu^T$. Since $u^Tu=1$, we have $H^T=H$ and
\[
H^2=(I_n-2uu^T)^2=I_n-4uu^T+4u(u^Tu)u^T=I_n.
\]
Thus $H^TH=I_n$, so $H$ is orthogonal. Moreover,
\[
D-cH=D-cI_n+2cuu^T=\diag(a_1,\ldots,a_n)+2cuu^T.
\]
We choose $u$ according to the parity of $n$.

\medskip
\noindent\textbf{Case 1: $n=2m+1$.}
Here $m\geq1$. Let
\[
\Delta=d_1-d_2+d_3-d_4+\cdots-d_{2m}+d_{2m+1}
\]
and let $L=c+\Delta$. We have $\Delta\geq0$ because
\[
\Delta=(d_1-d_2)+(d_3-d_4)+\cdots+(d_{2m-1}-d_{2m})+d_{2m+1}.
\]
Hence $L>0$. Let $x_j=a_{2j-1}$ for $j=1,\ldots,m+1$, and let $y_j=a_{2j}$ for $j=1,\ldots,m$. The inequalities among the $a_i$ give
\[
x_1>y_1>x_2>y_2>\cdots>x_m>y_m>x_{m+1}.
\]
Also $L>x_1$, since $L>0$ and $x_1=a_1<0$. Using $a_i=d_i-c$ and $L=c+\Delta$, we get
\[
L+\sum_{j=1}^{m}y_j-\sum_{j=1}^{m+1}x_j=c+\Delta+\sum_{j=1}^{m}(d_{2j}-c)-\sum_{j=1}^{m+1}(d_{2j-1}-c)=2c.
\]
Thus $L=x_1+\cdots+x_{m+1}+2c-y_1-\cdots-y_m$. Lemma \ref{lem:rank-one}, with $\sigma=2c$, gives a unit vector $w=(w_1,\ldots,w_{m+1})^T\in\R^{m+1}$ such that the eigenvalues of $\diag(a_1,a_3,\ldots,a_{2m+1})+2cww^T$, counting multiplicities, are exactly $L,a_2,a_4,\ldots,a_{2m}$.

Let $u\in\R^{2m+1}$ be defined by $u_{2j-1}=w_j$ for $j=1,\ldots,m+1$, and $u_{2j}=0$ for $j=1,\ldots,m$. Then $u^Tu=w^Tw=1$. Let $H=I_{2m+1}-2uu^T$ and let
\[
S=D-cH=\diag(a_1,\ldots,a_{2m+1})+2cuu^T.
\]
For each even coordinate vector $e_{2j}$, the vector $u$ has zero at each $(2j)$-th component. Hence $u^Te_{2j}=0$, and $Se_{2j}=a_{2j}e_{2j}$ for $j=1,\ldots,m$. Thus $a_2,a_4,\ldots,a_{2m}$ are eigenvalues of $S$ coming from the even coordinate directions. On the coordinates $1,3,\ldots,2m+1$, the matrix of $S$ is $\diag(a_1,a_3,\ldots,a_{2m+1})+2cww^T$. The eigenvalues on these coordinates are $L,a_2,a_4,\ldots,a_{2m}$. Hence the eigenvalues of $S$, counting multiplicities, are exactly one eigenvalue equal to $L$ and, for each $j=1,\ldots,m$, two eigenvalues equal to $a_{2j}$.

The matrix $S$ is symmetric, so the singular values of $S$ are the absolute values of the eigenvalues listed above. Since $L>0$ and $a_{2j}<0$, the singular values of $S$ are one value $L$ and, for each $j=1,\ldots,m$, two values equal to $c-d_{2j}$.

Since $d_2>d_4>\cdots>d_{2m}$, we have $c-d_{2m}>\cdots>c-d_4>c-d_2$. Moreover,
\[
L-(c-d_{2m})=\Delta+d_{2m}\geq0.
\]
Thus the singular values of $S$ can be listed in the first form in Lemma \ref{lem:two-term}. Lemma \ref{lem:two-term} gives $S\in\mathcal A_{2m+1}(2)$. Since $D=cH+S$ and $H\in O(2m+1)$, we get $D\in\mathcal A_{2m+1}(3)$.

\medskip
\noindent\textbf{Case 2: $n=2m$.}
Here $m\geq2$. Let
\[
\Delta=d_2-d_3+d_4-d_5+\cdots-d_{2m-1}+d_{2m}
\]
and let $L=c+\Delta$. Then
\[
\Delta=(d_2-d_3)+(d_4-d_5)+\cdots+(d_{2m-2}-d_{2m-1})+d_{2m}\geq0,
\]
so $L>0$. Let $x_j=a_{2j}$ for $j=1,\ldots,m$, and let $y_j=a_{2j+1}$ for $j=1,\ldots,m-1$. Then
\[
x_1>y_1>x_2>y_2>\cdots>x_{m-1}>y_{m-1}>x_m,
\]
and $L>x_1$ because $L>0$ and $x_1=a_2<0$. As in Case 1, using $a_i=d_i-c$ and $L=c+\Delta$, we get
\[
L+\sum_{j=1}^{m-1}y_j-\sum_{j=1}^{m}x_j=2c.
\]
Thus $L=x_1+\cdots+x_m+2c-y_1-\cdots-y_{m-1}$. Lemma \ref{lem:rank-one}, with $\sigma=2c$, gives a unit vector $w=(w_1,\ldots,w_m)^T\in\R^m$ such that the eigenvalues of $\diag(a_2,a_4,\ldots,a_{2m})+2cww^T$, counting multiplicities, are exactly $L,a_3,a_5,\ldots,a_{2m-1}$.

Let $u\in\R^{2m}$ be defined by $u_{2j}=w_j$ and $u_{2j-1}=0$ for $j=1,\ldots,m$. Let $H=I_{2m}-2uu^T$ and let
\[
S=D-cH=\diag(a_1,\ldots,a_{2m})+2cuu^T.
\]
The argument from Case 1, with the even and odd coordinates interchanged, gives the following eigenvalues of $S$. The odd coordinate directions give $a_1,a_3,\ldots,a_{2m-1}$, while the even coordinate block gives $L,a_3,a_5,\ldots,a_{2m-1}$. Therefore the eigenvalues of $S$, counting multiplicities, are exactly one eigenvalue equal to $L$, one eigenvalue equal to $a_1$, and, for each $j=1,\ldots,m-1$, two eigenvalues equal to $a_{2j+1}$.

Since $S$ is symmetric, the singular values of $S$ are the absolute values of these eigenvalues. Since $L>0$ and all $a_i<0$, the singular values are one value $L$, one value $c-d_1$, and, for each $j=1,\ldots,m-1$, two values equal to $c-d_{2j+1}$.

Since $d_1>d_3>d_5>\cdots>d_{2m-1}$, we have $c-d_{2m-1}>\cdots>c-d_5>c-d_3>c-d_1$. Also
\[
L-(c-d_{2m-1})=\Delta+d_{2m-1}\geq0.
\]
Thus the singular values of $S$ can be listed in the second form in Lemma \ref{lem:two-term}. Lemma \ref{lem:two-term} gives $S\in\mathcal A_{2m}(2)$. Since $D=cH+S$ and $H\in O(2m)$, we get $D\in\mathcal A_{2m}(3)$.
\end{proof}

\medskip

\begin{lemma}\label{lem:nonstrict-diagonal}
Let $n\geq3$, and let
\[
D=\diag(d_1,\ldots,d_n),\qquad d_1\geq d_2\geq\cdots\geq d_n\geq0.
\]
Then $D\in\mathcal A_n(3)$.
\end{lemma}

\begin{proof}
For each positive integer $N$, let
\[
d_i^{(N)}=d_i+\frac{n-i+1}{N+n}\qquad (i=1,\ldots,n)
\]
and let $D_N=\diag(d_1^{(N)},\ldots,d_n^{(N)})$. Then
\[
d_i^{(N)}-d_{i+1}^{(N)}=(d_i-d_{i+1})+\frac{1}{N+n}>0
\]
for $i=1,\ldots,n-1$. Hence $d_1^{(N)}>d_2^{(N)}>\cdots>d_n^{(N)}\geq0$. Also $D_N\to D$ as $N\to\infty$.

Let $c=d_1+2$. Since $d_1^{(N)}<d_1+1$, we have $c>d_1^{(N)}$ for every $N$. Applying Lemma \ref{lem:strict-diagonal} to $D_N$ with this same value of $c$, we obtain a unit vector $u_N\in\R^n$, a Householder reflection
\[
H_N=I_n-2u_Nu_N^T,
\]
and a matrix $S_N=D_N-cH_N$ whose singular values satisfy one of the two hypotheses in Lemma \ref{lem:two-term}.

Let $R_N=s_1(S_N)$ and $r_N=s_n(S_N)$, where $s_1(S_N)$ means the largest singular value of $S_N$, and $s_n(S_N)$ means the smallest singular value of $S_N$. By the final assertion of Lemma \ref{lem:two-term}, there are orthogonal matrices $U_N,V_N\in O(n)$ such that
\[
S_N=\frac{R_N+r_N}{2}U_N+\frac{R_N-r_N}{2}V_N.
\]
Let
\[
\alpha_N=\frac{R_N+r_N}{2},\qquad \beta_N=\frac{R_N-r_N}{2}.
\]
Then
\[
D_N=cH_N+\alpha_NU_N+\beta_NV_N.
\]

The sequences $\alpha_N$ and $\beta_N$ are bounded. This is because the spectral norm $\|\cdot\|_2$ is the largest singular value, so $R_N=\|S_N\|_2$. From $S_N=D_N-cH_N$, the triangle inequality gives
\[
R_N\leq \|D_N\|_2+c\|H_N\|_2.
\]
The matrix $H_N$ is orthogonal, so $\|H_N\|_2=1$. The matrix $D_N$ is nonnegative diagonal with decreasing diagonal entries, so $\|D_N\|_2=d_1^{(N)}<d_1+1$. Hence $R_N<d_1+1+c$. Since $0\leq r_N\leq R_N$, we have $0\leq\alpha_N\leq R_N$ and $0\leq\beta_N\leq R_N$. Therefore $\alpha_N$ and $\beta_N$ are bounded.

The unit sphere in $\R^n$ is compact. After passing to a subsequence, we may assume that $u_N\to u$ for some unit vector $u\in\R^n$. Hence
\[
H_N=I_n-2u_Nu_N^T\to I_n-2uu^T.
\]
Let $H=I_n-2uu^T$. Thus the limiting matrix $H$ is also a Householder reflection.

The group $O(n)$ is compact. After passing to a further subsequence, we may assume that $U_N\to U$ and $V_N\to V$ for some $U,V\in O(n)$. Since $\alpha_N$ and $\beta_N$ are bounded, after passing to a further subsequence we may assume that $\alpha_N\to\alpha$ and $\beta_N\to\beta$ for some real numbers $\alpha,\beta$.

Letting $N\to\infty$ in
\[
D_N=cH_N+\alpha_NU_N+\beta_NV_N
\]
gives
\[
D=cH+\alpha U+\beta V.
\]
The matrices $H,U,V$ are orthogonal. Hence $D\in\mathcal A_n(3)$.
\end{proof}

\medskip

\begin{theorem}\label{thm:main}
For every integer $n\geq3$,
\[
M_n(\R)=\mathcal A_n(3).
\]
Equivalently, every real $n\times n$ matrix with $n\geq3$ is a real linear combination of three real orthogonal matrices.
\end{theorem}

\begin{proof}
Let $A\in M_n(\R)$. By the real singular value decomposition, there are orthogonal matrices $P,Q\in O(n)$ and a diagonal matrix
\[
D=\diag(d_1,\ldots,d_n),\qquad d_1\geq\cdots\geq d_n\geq0,
\]
such that $A=PDQ^T$. By Lemma \ref{lem:nonstrict-diagonal}, there are real numbers $\gamma_1,\gamma_2,\gamma_3$ and orthogonal matrices $R_1,R_2,R_3\in O(n)$ such that
\[
D=\gamma_1R_1+\gamma_2R_2+\gamma_3R_3.
\]
Multiplying this identity on the left by $P$ and on the right by $Q^T$ gives
\[
A=\gamma_1(PR_1Q^T)+\gamma_2(PR_2Q^T)+\gamma_3(PR_3Q^T).
\]
For each $j=1,2,3$, the matrix $PR_jQ^T$ is orthogonal because
\[
(PR_jQ^T)^T(PR_jQ^T)=QR_j^TP^TPR_jQ^T=QR_j^TR_jQ^T=QQ^T=I_n.
\]
Therefore $A\in\mathcal A_n(3)$.
\end{proof}

\medskip

\begin{corollary}\label{cor:classification}
Let $k_{\min}(n)$ denote the least number of orthogonal summands needed to represent every real $n\times n$ matrix as a real linear combination of orthogonal matrices. Then the exact values are
\[
k_{\min}(1)=1,
\qquad
k_{\min}(2)=2,
\qquad
k_{\min}(n)=3\quad (n\geq3).
\]
\end{corollary}

\begin{proof}
When $n=1$, every real $1\times1$ matrix has the form $[a]$ with $a\in\R$, and $[a]=a[1]$. Since $[1]\in O(1)$, we have $k_{\min}(1)=1$.

Now let $n=2$. Let $A\in M_2(\R)$. By the real singular value decomposition, there are $P,Q\in O(2)$ and numbers $s_1\geq s_2\geq0$ such that
\[
A=P\diag(s_1,s_2)Q^T.
\]
The diagonal matrix satisfies
\[
\diag(s_1,s_2)=\frac{s_1+s_2}{2}I_2+\frac{s_1-s_2}{2}\diag(1,-1).
\]
Both $I_2$ and $\diag(1,-1)$ are orthogonal. Multiplying on the left by $P$ and on the right by $Q^T$, we get a representation of $A$ as a real linear combination of two orthogonal matrices. Therefore $M_2(\R)=\mathcal A_2(2)$.

One term does not suffice in dimension two. This is because if $A=\alpha Q$ with $Q\in O(2)$, then both singular values of $A$ are equal to $|\alpha|$. The matrix $\diag(1,0)$ has singular values $1$ and $0$, so it is not a scalar multiple of an orthogonal matrix. Hence $k_{\min}(2)=2$.

For $n\geq3$, Theorem \ref{thm:main} gives $k_{\min}(n)\leq3$. The lower bound $k_{\min}(n)>2$ is due to Li, Zhang, and Li \cite[Proposition~4.10]{LiZhangLi2022}. Therefore $k_{\min}(n)=3$ for every $n\geq3$.
\end{proof}

\medskip

\appendix
\section{Proofs of the preliminary lemmas}

\begin{proof}[Proof of Lemma \ref{lem:two-term}]
We first record a two-dimensional calculation. Let $a\geq b\geq0$, and let
\[
C(\theta)=
\begin{pmatrix}
\cos\theta&-\sin\theta\\
\sin\theta&\cos\theta
\end{pmatrix}.
\]
Then $C(\theta)$ is orthogonal. Let $B_\theta=aI_2+bC(\theta)$. If $x=a+b\cos\theta$ and $y=b\sin\theta$, then
\[
B_\theta=\begin{pmatrix}x&-y\\ y&x\end{pmatrix}.
\]
A direct multiplication gives
\[
B_\theta^TB_\theta=(x^2+y^2)I_2.
\]
Also $x^2+y^2=a^2+b^2+2ab\cos\theta$. Thus both singular values of $B_\theta$ are equal to
\[
(a^2+b^2+2ab\cos\theta)^{1/2}.
\]

Assume first that $n=2m+1$ with $m\geq1$ and that the singular values of $X$ can be listed as
\[
R,t_1,t_1,\ldots,t_m,t_m,
\]
where $R\geq t_1\geq\cdots\geq t_m\geq0$. Let $r=t_m$. Let
\[
a=\frac{R+r}{2},\qquad b=\frac{R-r}{2}.
\]
Then $a+b=R$ and $a-b=r$.

For each $j=1,\ldots,m$, we choose an angle $\theta_j$. If $b=0$, then $R=r$ and the inequalities force $t_j=R$ for all $j$; let $\theta_j=0$. If $b>0$ and $t_j=R$, let $\theta_j=0$. If $b>0$ and $t_j=r$, let $\theta_j=\pi$. If $b>0$ and $r<t_j<R$, then
\[
-1<\frac{t_j^2-a^2-b^2}{2ab}<1,
\]
so there is an angle $\theta_j\in(0,\pi)$ such that
\[
\cos\theta_j=\frac{t_j^2-a^2-b^2}{2ab}.
\]
In all cases, the identity $t_j^2=a^2+b^2+2ab\cos\theta_j$ holds. By the two-dimensional calculation above, the block $aI_2+bC(\theta_j)$ has singular values $t_j,t_j$.

Let
\[
W=\diag(1,C(\theta_1),\ldots,C(\theta_m)).
\]
Then $W\in O(2m+1)$. Let $B=aI_{2m+1}+bW$. The first block of $B$ gives the singular value $a+b=R$. The $j$th two-dimensional block gives the singular values $t_j,t_j$. Hence the singular values of $B$, counted with multiplicity, are exactly the same as those of $X$.

Now assume that $n=2m$ with $m\geq2$ and that the singular values of $X$ can be listed as
\[
R,t_1,t_1,\ldots,t_{m-1},t_{m-1},r,
\]
where $R\geq t_1\geq\cdots\geq t_{m-1}\geq r\geq0$. Let
\[
a=\frac{R+r}{2},\qquad b=\frac{R-r}{2}.
\]
Then $a+b=R$ and $a-b=r$.

For each $j=1,\ldots,m-1$, we choose an angle $\theta_j$. If $b=0$, then $R=r$ and the inequalities force $t_j=R$ for all $j$; let $\theta_j=0$. If $b>0$ and $t_j=R$, let $\theta_j=0$. If $b>0$ and $t_j=r$, let $\theta_j=\pi$. If $b>0$ and $r<t_j<R$, then
\[
-1<\frac{t_j^2-a^2-b^2}{2ab}<1,
\]
so there is an angle $\theta_j\in(0,\pi)$ such that
\[
\cos\theta_j=\frac{t_j^2-a^2-b^2}{2ab}.
\]
In all cases, the identity $t_j^2=a^2+b^2+2ab\cos\theta_j$ holds. Again the block $aI_2+bC(\theta_j)$ has singular values $t_j,t_j$.

Let
\[
W=\diag(1,C(\theta_1),\ldots,C(\theta_{m-1}),-1).
\]
Then $W\in O(2m)$. Let $B=aI_{2m}+bW$. The first block gives the singular value $a+b=R$, the last block gives the singular value $a-b=r$, and the two-dimensional blocks give the paired singular values $t_j,t_j$. Hence the singular values of $B$, counted with multiplicity, are exactly the same as those of $X$.

In either parity, we have constructed a matrix $B=aI_n+bW$ with $W\in O(n)$ and with the same singular values as $X$. Let
\[
X=P_X\Sigma Q_X^T,
\qquad
B=P_B\Sigma Q_B^T
\]
be real singular value decompositions with the same nonnegative diagonal matrix $\Sigma$. Then
\[
X=(P_XP_B^T)B(Q_BQ_X^T).
\]
Using $B=aI_n+bW$, we get
\[
X=a(P_XP_B^TQ_BQ_X^T)+b(P_XP_B^TWQ_BQ_X^T).
\]
Both matrices multiplying $a$ and $b$ are products of orthogonal matrices. Hence both are orthogonal. Thus $X\in\mathcal A_n(2)$. Since $a=(R+r_X)/2$ and $b=(R-r_X)/2$, where $r_X=s_n(X)$, the final assertion also follows.
\end{proof}

\medskip

\begin{proof}[Proof of Lemma \ref{lem:rank-one}]
Let
\[
p(\lambda)=(\lambda-L)\prod_{j=1}^{\ell-1}(\lambda-y_j).
\]
For each $i=1,\ldots,\ell$, let
\[
w_i=\left(-\frac{p(x_i)}{\sigma\prod_{k\neq i}(x_i-x_k)}\right)^{1/2}.
\]

We check that the number under the square root is positive. Fix $i$. Since $L>x_1\geq x_i$, the factor $x_i-L$ is negative. Among the numbers $y_1,\ldots,y_{\ell-1}$, exactly $y_1,\ldots,y_{i-1}$ are larger than $x_i$, and the remaining $y_i,\ldots,y_{\ell-1}$ are smaller than $x_i$. Hence $\prod_{j=1}^{\ell-1}(x_i-y_j)$ has sign $(-1)^{i-1}$. Therefore $p(x_i)$ has sign $(-1)^i$. Similarly, $\prod_{k\neq i}(x_i-x_k)$ has sign $(-1)^{i-1}$, since exactly the $i-1$ factors with $k<i$ are negative. Since $\sigma>0$, the quotient
\[
-\frac{p(x_i)}{\sigma\prod_{k\neq i}(x_i-x_k)}
\]
is positive. Thus $w_i$ is real and nonzero.

Let $X=\diag(x_1,\ldots,x_\ell)$ and let $w=(w_1,\ldots,w_\ell)^T$. Define
\[
q(\lambda)=\det(\lambda I_\ell-X-\sigma ww^T).
\]
For $\lambda\notin\{x_1,\ldots,x_\ell\}$, the matrix determinant lemma (see Horn and Johnson \cite{HornJohnson2013}) gives
\[
q(\lambda)=\prod_{i=1}^\ell(\lambda-x_i)-\sigma\sum_{i=1}^\ell w_i^2\prod_{k\neq i}(\lambda-x_k).
\]
Both sides are polynomials in $\lambda$, so this identity extends to every real $\lambda$. At $\lambda=x_i$, every term in the sum vanishes except the $i$th term. Therefore
\[
q(x_i)=-\sigma w_i^2\prod_{k\neq i}(x_i-x_k)=p(x_i).
\]
The polynomials $q(\lambda)$ and $p(\lambda)$ are both monic of degree $\ell$. Thus $q(\lambda)-p(\lambda)$ has degree at most $\ell-1$. It has the $\ell$ distinct zeros $x_1,\ldots,x_\ell$, so it is the zero polynomial. Hence $q(\lambda)=p(\lambda)$.

The coefficient of $\lambda^{\ell-1}$ in $q(\lambda)$ is $-(x_1+\cdots+x_\ell)-\sigma(w_1^2+\cdots+w_\ell^2)$. The coefficient of $\lambda^{\ell-1}$ in $p(\lambda)$ is $-(L+y_1+\cdots+y_{\ell-1})$. Since $L+y_1+\cdots+y_{\ell-1}=x_1+\cdots+x_\ell+\sigma$, comparison of these coefficients gives $w_1^2+\cdots+w_\ell^2=1$. Thus $w$ is a unit vector.

Finally, since $q(\lambda)=p(\lambda)$, the characteristic polynomial of $X+\sigma ww^T$ is \[(\lambda-L)\prod_{j=1}^{\ell-1}(\lambda-y_j).\] Therefore the eigenvalues of $X+\sigma ww^T$, counting multiplicities, are exactly $L,y_1,\ldots,y_{\ell-1}$.
\end{proof}


\begin{thebibliography}{9}

\bibitem{BunchNielsenSorensen1978}
J.~R. Bunch, C.~P. Nielsen, and D.~C. Sorensen.
Rank-one modification of the symmetric eigenproblem.
\emph{Numerische Mathematik}, 31:31--48, 1978/79.

\bibitem{Golub1973}
G.~H. Golub.
Some modified matrix eigenvalue problems.
\emph{SIAM Review}, 15(2):318--334, 1973.

\bibitem{HornJohnson2013}
R.~A. Horn and C.~R. Johnson.
\emph{Matrix Analysis}.
Cambridge University Press, Cambridge, second edition, 2013.

\bibitem{LiPoon2002}
C.-K. Li and E. Poon.
Additive decomposition of real matrices.
\emph{Linear and Multilinear Algebra}, 50(4):321--326, 2002.

\bibitem{LiZhangLi2022}
Z.~Li, T.~Zhang, and C.-J. Li.
Singular value structure of real matrices which can be expressed as a linear combination of two orthogonal matrices.
\emph{Operators and Matrices}, 16(3):775--792, 2022.

\bibitem{Wu1994}
P.~Y. Wu.
Additive combinations of special operators.
\emph{Banach Center Publications}, 30(1):337--361, 1994.

\bibitem{Zhan2001}
X.~Zhan.
Span of the orthogonal orbit of real matrices.
\emph{Linear and Multilinear Algebra}, 49(4):337--346, 2001.

\bibitem{Zhan2007}
X.~Zhan.
Open problems in matrix theory.
In \emph{Proceedings of the International Congress of Chinese Mathematicians 2007}, Vol.~II, pages 1--17, 2007.

\end{thebibliography}
\end{document}